\documentclass{elsarticle}
 \usepackage{graphicx,amssymb,amstext,amsmath}
 \usepackage{bigstrut}
 \usepackage{enumerate}
 \usepackage{amsthm}
 \usepackage[english]{babel}
 \usepackage[all]{xy}
 \usepackage[vcentermath]{youngtab}
\usepackage{longtable}
\usepackage{lscape}

 \newtheorem{theorem}{Theorem}[section]
 \newtheorem{thm}[theorem]{Theorem}
 
 \newtheorem{lem}[theorem]{Lemma}
 \newtheorem{cor}[theorem]{Corollary}
 \newtheorem{prop}[theorem]{Proposition}

\newcommand{\zz}{\mathbb{Z}}

\DeclareMathOperator\id{id}
\DeclareMathOperator\fix{fix}

\DeclareMathOperator\alt{Alt}
\DeclareMathOperator\sym{Sym}
\DeclareMathOperator\PGL{PGL}
\DeclareMathOperator\PSL{PSL}

\DeclareMathOperator\cc{\mathcal{C}}
\DeclareMathOperator\dd{\mathcal{D}}
\DeclareMathOperator\one{\bf{1}}

\begin{document}
\title{The Erd\H{o}s-Ko-Rado property for some $2$-transitive groups}

 \author{Bahman Ahmadi}
 \ead{ahmadi2b@uregina.ca}

 \author{Karen Meagher \corref{cor1}\fnref{fn1}}
 \ead{karen.meagher@uregina.ca}
 \cortext[cor1]{Corresponding author}
 \fntext[fn1]{Research supported by NSERC.}

 \address{Department of Mathematics and Statistics,\\
 University of Regina, 3737 Wascana Parkway, S4S 0A4 Regina SK, Canada}

 \begin{abstract}
   A subset of a group $G \leq \sym(n)$ is \textsl{intersecting} if
   for any pair of permutations $\pi,\sigma \in G$ there is an $i \in
   \{1,2,\dots,n\}$ such that $\pi(i) = \sigma(i)$. It has been shown,
   using an algebraic approach, that the largest intersecting sets in
   each of $\sym(n)$, $\alt(n)$ and $\PGL(2,q)$ are exactly the cosets
   of the point-stabilizers. In this paper, we show how this method
   can be applied more generally to many $2$-transitive groups. We
   then apply this method to the Mathieu groups and to all 2-transtive
   groups with degree no more than 20.
 \end{abstract}


\maketitle

\section{Introduction}\label{introduction}

There have been several recent publications~\cite{2013arXiv1302.7313A,
  MR2009400, Karen, KuW07, MR2061391, MR2876903, MeagherS11} that
determine the maximum sets of elements from a permutation group such
that any two permutations from the set both map at least one point to
the same element.  These results are considered to be versions of the
Erd\H{o}s-Ko-Rado (EKR) theorem~\cite{TheOriginalEKR} for
permutations.  In this paper, we focus on generalizing the algebraic
method used in \cite{Karen, MeagherS11} to prove the natural version
the EKR theorem for $2$-transitive permutation groups.

The EKR theorem describes the size and structure of the largest
collection of intersecting $k$-subsets of the set
$[n]=\{1,2,\dots,n\}$. Provided that $n >2k$, the largest such collection
is comprised of all $k$-subsets that contain a common fixed element
from $[n]$. This result has been extended to other objects
for which there is a notion of intersection. For examples, an
analogous result has been shown to hold for integer
sequences~\cite{MR1722210}, vector spaces over a finite
field~\cite{MR0382015}, partitions~\cite{MR2156694} and many other
objects. Here we consider the extension to permutations.

Let $G\leq \sym(n)$ be a permutation group with the natural action on
the set $[n]$. Two permutations $\pi,\sigma\in G$ are said to
\textsl{intersect} if $\pi\sigma^{-1}$ has a fixed point in $[n]$. In
other words, $\sigma$ and $\pi$ do not intersect if $\pi\sigma^{-1}$
is a \textsl{derangement}.  Throughout this paper, we denote the set
of all derangement elements of group the $G$ by $\dd_G$.  A subset
$S\subseteq G$ is, called \textsl{intersecting} if any pair of its
elements intersect. Clearly, the stabilizer of a point is an
intersecting set in $G$; as is any coset of the stabilizer of a
point. Note that such a coset is the collection of all permutations in
the group that map $i$ to $j$ for some $i,j \in [n]$. In this paper we
focus on $2$-transitive groups and if $G$ is any $2$-transitive subgroup of
$\sym(n)$, then any stabilizer of a point has size $\frac{|G|}{n}$.

We say the group $G$ has the \textsl{EKR property}, if the size of any
intersecting subset of $G$ is bounded above by the size of the largest
point-stabilizer in $G$. Further, $G$ is said to have the
\textsl{strict EKR property} if the only maximum intersecting subsets
of $G$ are the cosets of the point-stabilizers. It is clear from the
definition that if a group has the strict EKR property, then it will
also have the EKR property.

It has been shown that the symmetric group has the strict EKR
property~\cite{MR2009400, Karen, MR2061391}, as does the alternating
group~\cite{2013arXiv1302.7313A,KuW07} and the group
$\PGL(2,q)$~\cite{MeagherS11}. In this paper we focus on an the algebraic
method that was used in~\cite{Karen, 2013arXiv1302.7313A,
  MeagherS11} to show that these groups have the strict EKR
property. First we show how this method is particularly effective for
$2$-transitive groups. In fact, we can prove a group has the strict
EKR property using a few straight-forward calculations. As an example
of this, we show that all the $2$-transitive Mathieu groups have this
property. We also consider all the $2$-transtive groups with degree no
more than 20; all of these have the EKR property and we show that
several also have the strict EKR property.

\section{Derangement Graphs}

For a group $G$ we identify a graph such that the independent sets in
this graph are intersecting sets of permutations from the group. This
graph is called the \textsl{derangement graph} for $G$. The vertices
of this graph are the elements in $G$ and two vertices are adjacent if
and only if they are not intersecting. This graph is the Cayley graph
on $G$ with the derangements from $G$, so the set $\dd_G$, as the
connection set. We denote the derangement graph by $\Gamma_G
=\Gamma(G,\dd_G)$.  The degree of $\Gamma_G$ is $|\dd_G|$ and it is a
vertex-transitive graph.

For any $i,j\in [n]$, we define the \textsl{canonical independent
  sets} to be the sets
\[
S_{i,j}=\{\pi\in G\,\,|\,\, \pi(i)=j\}.
\]
These are clearly independent sets in $\Gamma_G$ and they are the
cosets of the point-stabilizers in $G$ under the natural action of $G$
on $[n]$. For each $i,j\in [n]$ it is easy to see that
$|S_{i,j}|=\frac{|G|}{n}$. We will denote the characteristic vector of
$S_{i,j}$ by $v_{i,j}$.

The group $G$ has the EKR property if 
\[
\alpha(\Gamma_G) = \frac{|G|}{n}.
\]
Furthermore, a group has the strict EKR property if every independent
set of maximum size is equal to $S_{i,j}$ for some $i$ and $j$.

For any group there is an association scheme called the
\textsl{conjugacy class scheme}.  The matrices $\{A_1,\dots,
A_k\}$ in the conjugacy class scheme of $G$ are the $|G|\times |G|$
matrices such that, for any $1\leq i\leq k$, the entry $(g,h)$ of
$A_i$ is $1$ if $hg^{-1}$ belongs to the $i$-th conjugacy class,
and $0$ otherwise. The derangement graph for any group $G$ is a union
of graphs in the conjugacy class scheme of $G$. (For more details on
association schemes and the conjugacy class scheme see~\cite[Example
2.1 (2)]{MR882540}.)

The derangement graph is a normal Cayley graph (meaning that the
connection set is closed under conjugation).  By the eigenvalues,
eigenvectors and eigenspaces of a graph we mean the eigenvalues,
eigenvectors and eigenspaces of the adjacency matrix of the graph.
Using a result of Diaconis and Shahshahani~\cite{MR626813}, it is
possible to calculate the eigenvalues of the derangement graph using
the irreducible characters of the group (see \cite[Chapter
4]{BahmanAhmadi} for a detailed proof of this).

\begin{cor}\label{evals_of_der_graph} 
The eigenvalues of the derangement graph $\Gamma_{G}$ are given by
\[
\eta_{\chi}=\frac{1}{\chi(\id)} \sum_{x \in \dd_G} \chi(x), 
\]
where $\chi$ ranges over all the distinct irreducible
characters of $G$. The multiplicity of the eigenvalue $\eta_\chi$ is
equal to $\chi(\id)^2$.\qed
\end{cor}

For any group $G$ the largest eigenvalue of $\Gamma_G$ is $|\dd_G|$
and the all ones vector is an eigenvector; this is the eigenvalue
given by the trivial character.

In this paper we focus on the \textsl{standard character}.  The value
of the standard character evaluated on any permutation is one less
than the number of elements fixed by that permutation. This is indeed
a character for any group, but for the $2$-transitive groups it is an
irreducible character (this follows by Burnside's lemma,
see~\cite[Chapter 3]{BahmanAhmadi} for a detailed proof).

\begin{lem}\label{standevalue}
  Let $G$ be a $2$-transitive subgroup of $\sym(n)$, then
  $-\frac{|\dd_G|}{n-1}$ is an eigenvalue of $\Gamma_G$.
\end{lem}
\proof By definition, the value of the standard character on the
identity of $G$ is $n-1$ and its value on any derangement is $-1$.  Applying
Corollary~\ref{evals_of_der_graph}, we have that the eigenvalue for the standard
character, $\chi$, is
\[
\eta_\chi = \frac{1}{\chi(\id)} \sum_{g \in \dd_G} \chi(g) = -\frac{|\dd_G|}{n-1}.\qed
\]

This implies that the multiplicity of $-\frac{|\dd_G|}{n-1}$, as an
eigenvalue, is at least $(n-1)^2$. If another irreducible character
also gives this eigenvalue, then the multiplicity will be larger.

\section{EKR property for $2$-transitive groups}

We focus on two methods that have been used to prove that a group has
the EKR property. The first method uses a bound called the
\textsl{clique-coclique bound} (\textsl{coclique} is another term for an
independent set in a graph). The version we use here was originally
proved by Delsarte \cite{MR0384310}. Before stating this bound, we
will need some notation.

Assume $\mathcal{A}=\{A_0, A_1,\ldots, A_d\}$ is an association scheme
on $v$ vertices and let $\{E_0, E_1,\ldots, E_d\}$ be the idempotents
of this association scheme. (For a detailed discussion about
association schemes, the reader may refer to \cite{MR2047311} or
\cite{MR882540}.)  For our purposes here, we only need to know that
the matrices $A_i$ are simultaneously diagonalizable, and the $E_j$
are the projections to the common eigenspaces of the matrices in the
association scheme. A \textsl{graph in an association scheme} $\mathcal{A}$ is a
graph whose adjacency matrix is one of the matrices in $\mathcal{A}$.

For a set $S$ of vertices in a graph $X$, the characteristic vector of
$S$ will be denoted by $v_S$ (the entries of $v_S$ are indexed by
$V(X)$ and the $v$-entry is equal to $1$ if $v \in S$ and $0$
otherwise).

\begin{thm}\label{clique_coclique_bound} 
  Let $X$ be the union of some of the graphs in an association scheme
  $\mathcal{A}$ on $v$ vertices. If $C$ is a clique and $S$ is an
  independent set in $X$, then
\[
|C||S|\leq v.
\]
If equality holds then
\[
v_C^\top\,E_j\,v_C\,\,v_S^\top\,E_j\,v_S = 0,\quad \text{for all} \quad j > 0.\qed
\]
\end{thm}

This gives a simple way to check if a group has the EKR property.

\begin{cor}
If a $2$-transitive group has a sharply-transitive set, then it has the EKR property.
\end{cor}
\proof A sharply-transitive set in a group $G \leq \sym(n)$ is a
clique of size $n$ in $\Gamma_G$. By the ratio bound, the size of the
largest independent set in $\Gamma_G$ is $\frac{|G|}{n}$. \qed

The clique-coclique bound actually holds for any vertex-transitive
graph, but we are interested in the above version for association
schemes because of the two following simple, but useful, corollaries
which were proven in~\cite{Karen}.

\begin{cor}\label{clique_vs_coclique}
  Let $X$ be a union of graphs in an association scheme such that the
  clique-coclique bound holds with equality in $X$.  Assume that $C$
  is a maximum clique and $S$ is a maximum independent set in
  $X$. Then, for $j > 0$, at most one of the vectors $E_jv_C$ and
  $E_jv_S$ is not zero.\qed
\end{cor}

In other words, provided that the clique-coclique bound holds with
equality, for any module of $\Gamma_G$ (other than the trivial module)
the projection of at most one of the vectors $v_C$ and $v_S$ will be
non-zero, where $S$ is any maximum independent set and $C$ is any
maximum clique.

For any group $G$ the clique-coclique bound
(Theorem~\ref{clique_coclique_bound}) applies to $\Gamma_{G}$ as
$\Gamma_G$ is the union of graphs in the conjugacy class scheme. The
idempotents of this scheme are $E_\chi$ where
\begin{equation}\label{idempotent}
(E_\chi)_{\pi,\sigma}=\frac{\chi(1)}{|G|}\chi(\pi^{-1}\sigma).
\end{equation}
where $\chi$ runs through the set of all irreducible characters of
$G$.  The vector space generated by the columns of $E_\chi$ is called
the \textsl{$\chi$-module} of $\Gamma_G$.

For any character $\chi$ of $G$ and any subset $X$ of $G$ define
\[
\chi(X)=\sum_{x\in X}\chi(x).
\]
Using Corollary~\ref{clique_vs_coclique} and Equation (\ref{idempotent}) one observes the following.

\begin{cor}\label{at_most_one_non-zero}
  Assume the clique-coclique bound holds with equality for the graph
  $\Gamma_G$ and let $\chi$ be an irreducible character of $G$ that is
  not the trivial character. If there is a clique $C$ of maximum size
  in $\Gamma_G$ with $\chi(C)\neq 0$, then
\[
E_\chi\,v_S = 0
\]
for any maximum independent set $S$ of $\Gamma_G$.\qed
\end{cor}

Let $G \leq \sym(n)$ and assume that for every irreducible
representation $\chi$ of $G$, except the standard representation and
the trivial representation, we can find a clique $C$ in $\Gamma_G$ of
size $n$ such that $\chi(C) \neq 0$. Then the above corollary implies
that the charateristic vector of any maximum independent set is in the
span of the trivial module and the standard module.

The other method we use to show that a group has the EKR property is
an eigenvalue bound called the \textsl{ratio bound}. The ratio bound
is due to Delsarte who used a linear programming argument to prove
this if for association schemes (see \cite[ Section 3.2]{Newman}).

\begin{thm}\label{ratio2}
  Let $X$ be a $k$-regular graph on $n$ vertices with $\tau$ the least
  eigenvalue of $X$. For any independent set $S$ we have
\[
|S|\leq \frac{n}{1-\frac{k}{\tau}}.
\]
Furthermore, the equality holds if and only if
\[
A(X)\left(v_S-\frac{|S|}{n}\mathbf{1}\right)=\tau\left(v_S-\frac{|S|}{n}\mathbf{1}\right).\qed
\]
\end{thm}

There are large families of groups for which this bound can be used to
show that group has the EKR property.

\begin{lem}\label{standardleast}
  Let $G$ be a $2$-transitive group. If the eigenvalue arising from
  the standard representation of $G$ is the least eigenvalue of $\Gamma_G$, then $G$ has the EKR property.
\end{lem}
\proof
Simply apply the eigenvalue from Lemma~\ref{standevalue} in the ratio bound to get that 
\[
\alpha(\Gamma_G) \leq \frac{|G|}{1 - \frac{|\dd_G|}{\frac{|\dd_G|}{n-1}}} = \frac{|G|}{n}.\qed
\]

The eigenvalue arising from the standard representation is the least
eigenvalue for $\Gamma_G$ for the following groups:
$\sym(n)$~\cite{MR2365981}, $\PGL(2,q)$,
$\PSL(2,q)$~\cite{MeagherS11}, and all the $2$-transitive Mathieu
groups.  Further, in the appendix we have a list of all the
$2$-transitive groups with degree no more than 20 for which this
holds.

Note that, the second part of Theorem~\ref{ratio2} states that if
$-\frac{|\dd_G|}{n-1}$ is the least eigenvalue, then the
characteristic vector of any maximum independent set lies in the
direct sum of the $\dd_G$-eigenspace and the
$-\frac{|\dd_G|}{n-1}$-eigenspace of $\Gamma_G$. If the standard
representation is the only representation that gives the least
eigenvalue of $\Gamma_G$, then the characteristic vector of any
independent set, when shifted by the all ones vector, lies in the
standard module.


\section{Strict EKR theorem for $2$-transitive groups}

In this section we describe a method used to show that a
$2$-transitive group has the strict EKR property. We call this the
\textsl{module method}. This method has several components. First, the
group must have the EKR property. Second, the characteristic vector
for any maximum independent set must be in the sum of the standard module and
the trivial module. We will show that the vectors $v_{i,j}$ form a
spanning set for the sum of these two modules; hence the
characteristic vector of every maximum independent set is a linear
combination of this vectors. Finally, if the only linear combination
of these vectors that gives the characteristic vector of a maximum
independent set is one of the vectors $v_{i,j}$, then the strict EKR
theorem holds for the group.

Before stating the conditions we need to check to show a
$2$-transitive group has the strict EKR property, we will need some
technical lemmas. The first two give a subset of the vectors
$v_{i,j}$ form a basis for the standard module.

\begin{lem}\label{instandard}
  Let $G$ be a $2$-transitive group and $S_{i,j}$ the canonical
  independent sets of $G$. Then $v_{i,j} - \frac{1}{n} \one$ lies in the standard module. 
\end{lem}
\proof Let $\chi$ be the standard representation, we will show that
$E_\chi( v_{i,j} - \frac{1}{n} \one) = v_{i,j} - \frac{1}{n}
\one$. Since $G$ is $2$-transitive we can assume with out loss of
generality that $i=j=n$. First note that
\[
E_\chi (v_{n,n}-\frac{1}{n}\one) = E_\chi (v_{n,n}). 
\]
Denote the row of $E_\chi$ corresponding to $\pi$ by
$[E_\chi]_\pi$. If $\pi(n) = n$ then
\[
[E_\chi]_\pi \cdot v_{n,n} = \frac{n-1}{|G|}\sum_{\sigma(n) = n} \chi(\pi^{-1}\sigma) = \frac{n-1}{|G|}\sum_{g(n) = n} \chi(g).
\]
Use $\fix'$ to denote the number of fixed points of an element from $G_n$ has on $[n-1]$, 
then this is equal to 
\[
\frac{n-1}{|G|}\sum_{g \in G_n} \fix'(g) = \frac{n-1}{|G|}\left(\frac{|G|}{n}\right) = \frac{n-1}{n}.
\]
The first equality holds by Burnside's lemma and the fact that $G_n$ is transitive.

Note that since $\sum_{g \in G} \chi(g) = 0$ this implies that
\[
 \sum_{i=1}^{n-1}\sum_{g(n) = i} \chi(g) = -\frac{|G|}{n},
\]
and since $G$ is $2$-transitive, for any $i \neq n$
\[
\sum_{g(n) = i} \chi(g) = -\frac{|G|}{n(n-1)},
\]
We can apply this in the case that $\pi(n) \neq n$ to get that
\[
[E_\chi]_\pi \cdot v_{n,n} = \frac{n-1}{|G|}\sum_{\sigma(n) =i } \chi(\pi^{-1}\sigma) = \frac{n-1}{|G|}(-\frac{|G|}{n(n-1)}) = -\frac{1}{n}.\qed
\]

\begin{lem}\label{basis_for_standard}
  Let $G$ be a $2$-transitive group. The set
\[
B:=\{v_{i,j}-\frac{1}{n}\mathbf{1}\,|\, i,j\in[n-1]\} 
\]
is a basis for the standard module of $G$.
\end{lem}
\proof Let $V$ denote the standard module of $G$. According to
Lemma~\ref{instandard}, $B\subset V$ and since the dimension of $V$ is
equal to $|B|=(n-1)^2$, it suffices to show that $B$ is linearly
independent. Note, also, that since $\mathbf{1}$ is not in the span of
$v_{i,j}$ for $i,j\in[n-1]$, it is enough to prove that the set
$\{v_{i,j}\,|\, i,j\in[n-1]\}$ is linearly independent.

Define a matrix $L$ to have the vectors $v_{i,j}$, with $i,j\in[n-1]$,
as its columns.  Then the rows of $L$ are indexed by the elements of
$G$ and the columns are indexed by the ordered pairs $(i,j)$, where
$i,j\in [n-1]$; we will also assume that the ordered pairs are listed
in lexicographic order.  It is, then, easy to see that
\[
L^\top L=\frac{(n-1)!}{2}\,I_{(n-1)^2}\,+\, \frac{(n-2)!}{2}\left( A(K_{n-1})\otimes A(K_{n-1})\right),
\]
where $I_{(n-1)^2}$ is the $(n-1) \times (n-1)$ identity matrix,
$A(K_{n-1})$ is the adjacency matrix of the complete graph $K_{n-1}$
and $\otimes$ is the tensor product.  The distinct eigenvalues of
$A(K_{n-1})$ are $-1$ and $n-2$; thus the eigenvalues of
$A(K_{n-1})\otimes A(K_{n-1})$ are $-(n-2), 1, (n-2)^2$. This implies
that the least eigenvalue of $L^\top L$ is
\[
\frac{(n-1)!}{2}-\frac{(n-2)(n-2)!}{2}>0.
\]
This proves that $L^\top L$ is non-singular and hence full rank. This,
in turn, proves that $L$ has full rank and that $\{v_{i,j}\,|\,
i,j\in[n-1]\}$ is linearly independent.\qed

Define the $|G|\times n^2$ matrix $H$ to be the matrix whose columns
are the vectors $v_{i,j}$, for all $i,j\in [n]$. Note that since $H$
has constant row-sums, the vector $\mathbf{1}$ is in the column space
of $H$.  We denote by $H_{(i,j)}$ the column of $H$ indexed by the
pair $(i,j)$, for any $i,j\in [n]$. Define the matrix $\overline{H}$
to be the matrix obtained from $H$ by deleting all the columns
$H_{(i,n)}$ and $H_{(n,j)}$ for any $i,j\in[n-1]$. With a similar
method as in the proof of \cite[Proposition 10]{MeagherS11}, we prove
the following.
\begin{lem}\label{col_H_bar}
The matrices $H$ and $\overline{H}$ have the same column space.
\end{lem}
\proof Obviously, the column space of $\overline{H}$ is a subspace of
the column space of $H$; thus we only need to show that the vectors
$H_{(i,n)}$ and $H_{(n,j)}$ are in the column space of $\overline{H}$,
for any $i,j\in[n-1]$. Since $G$ is $2$-transitive, it suffices to
show this for $H_{(1,n)}$. Define the vectors $v$ and $w$ as follows:
\[
v:=\sum_{i\neq 1,n} \sum_{j\neq n} H_{(i,j)}\quad\text{and}\quad w:=(n-3) \sum_{j\neq n} H_{(1,j)}\,+H_{(n,n)}.
\]
The vectors $v$ and $w$ are in the column space of $\overline{H}$. It
is easy to see that for any $\pi \in G$,
\[
v_\pi=
\begin{cases}
n-2,& \quad\text{if}\quad \pi(1)=n;\\
n-2, & \quad \text{if} \quad\pi(n)=n;\\
n-3,& \quad\text{otherwise},
\end{cases}
\quad\quad\quad
w_\pi=
\begin{cases}
0,& \quad\text{if}\quad \pi(1)=n;\\
n-2, &  \quad\text{if}\quad\pi(n)=n;\\
n-3,& \quad\text{otherwise}.
\end{cases}
\]
Thus 
\[
(v-w)_\pi=
\begin{cases}
n-2,& \quad\text{if}\quad \pi(1)=n;\\
0, & \quad \text{if} \quad\pi(n)=n;\\
0,& \quad\text{otherwise},
\end{cases}
\]
which means that $(n-2)H_{(1,n)}=v-w$. This completes the proof.\qed

If the columns of $\overline{H}$ are arranged so that the first $n$
columns correspond to the pairs $(i,i)$, for $i\in [n]$, and the rows
are arranged so that the first row corresponds to the identity
element, and the next $|\dd_G|$ rows correspond to the derangements of
$G$, then $\overline{H}$ has the following block structure:
\[
\begin{bmatrix}
1& 0\\
0 & M \\
B & C\\
\end{bmatrix}.
\]
Note that the rows and columns of $M$ are indexed by the elements of
$\dd_G$ and the pairs $(i,j)$ with $i,j\in[n-1]$ and $i\neq j$,
respectively; thus $M$ is a $|\dd_G|\times (n-1)(n-2)$
matrix. Throughout the paper, we will refer to this matrix simply as
``the matrix $M$ for $G$''.

The next proposition shows that the submatrix $B$ in $\overline{H}$
above contains an $n \times n$ identity matrix.

\begin{prop}\label{condition_d} 
  Let $G\leq \sym(n)$ be $2$-transitive. Then for any $x\in [n]$,
  there is an element in $G$ which has $x$ as its only fixed point.
\end{prop}
\proof Since $G$ is transitive, it suffices to show this for $x=n$. We
need to show that the stabilizer of $n$ in $G$, denoted $G_n$, has a
derangement. Suppose for every element $g\in G_n$, we have
$|\fix(g)|\geq 1$. This means that
\[
\frac{1}{|G_n|} \sum_{g\in G_n} |\fix(g)|  \geq \frac{(n-1) +|G_n| -1}{|G_n|} = \frac{(n-2 + |G_n|)}{|G_n|},
\]
which is greater than $1$, if $n>2$. Hence by Burnside's lemma, the
number of orbits of the action of $G_n$ on $[n-1]$ is more
than one which is a contradiction since $G_n$ acts transitively on
$[n-1]$. Thus there must be a derangement in $G_n$.\qed

We are now ready to state the method that we call the \textsl{module method} which we use to prove
that many $2$-transitive groups have the strict EKR theorem.

\begin{thm}\label{module_method_thm} Let $G\leq \sym(n)$ be 2-transitive and assume the following conditions hold:
\begin{enumerate}[(a)]
\item $G$ has the EKR  property;
\item for any maximum intersecting set $S$ in $G$, the vector $v_S$ lies in the direct sum of the trivial and the standard modules of $G$; and 
\item the matrix $M$ for $G$ has full rank.
\end{enumerate}
Then $G$ has the strict EKR property.
\end{thm}
\proof Since $G$ has the EKR property, the maximum size of an
  intersecting subset of $G$ is $|G|/n$, i.e. the size of a
  point-stabilizer. Suppose that $S$ is an intersecting set of maximum size. It is enough
  to show that $S=S_{i,j}$, for some $i,j\in [n]$.  Without loss of
  generality, we may assume that $S$ includes the identity element.
  By assumption (b) and Lemma~\ref{basis_for_standard}, $v_S$ is
  in the column space of $H$; thus according to Lemma~\ref{col_H_bar},
  $v_S$ belongs to the column space of $\overline{H}$; therefore
\[
\begin{bmatrix}
1 & 0\\
0 & M \\
B & C\\
\end{bmatrix}\begin{bmatrix} z \\ w  \end{bmatrix}
=v_S
\]
for some vectors $z$ and $w$. Since the identity is in $S$, no elements from $\dd_G$ are in $S$, the characteristic vector of
$S$ has the form
\[
v_S= \begin{bmatrix} 1 \\ 0 \\ t  \end{bmatrix}
\]
for some vector $01$-vector $t$. Thus we have 
\[
1^\top z=1, \quad Mw=0, \quad Bz+Cw=t.
\]
Since $M$ has full rank, $w=0$ and so $Bz=t$. Furthermore,
according to Proposition~\ref{condition_d}, one can write
\[
B=\begin{bmatrix} I_{n} \\[.2cm] B'  \end{bmatrix}\quad 
\text{and}\quad 
Bz=\begin{bmatrix} z \\[.2cm] B'z  \end{bmatrix}. 
\]
Since $Bz$ is equal to the $01$-vector $t$, the vector $z$ must also
be a $01$-vector. But, on the other hand, $1^\top z=1$, thus we
conclude that exactly one of the entries of $z$ is equal to $1$. This
means that $v_S$ is the characteristic vector of the stabilizer of a
point.\qed

Through this paper we refer the conditions of
Theorem~\ref{module_method_thm} as conditions (a), (b) and (c),
without reference to the theorem.  We point out that the module
condition (b) is the reason we call this method the module method.

This gives an algorithm for testing if a $2$-transitive group has the
EKR or strict EKR property (although it cannot determine if a group
does not have the strict EKR property). First we calculate all the
eigenvalues of $\Gamma_G$. If the standard representation gives the
least eigenvalue, then the EKR property holds.  If it is the only
representation that gives the least eigenvalue, then the
characteristic vector for any maximum independent set is in the
standard module. If the matrix $M$ has full rank, then the group has
the strict EKR property. If we can't show that the standard
representation gives the least eigenvalue, then we check it the
derangement graph of the group has a clique of size $n$. If it does,
then the group has the EKR property. If for each irreducible character
(other than the trivial and the standard character) we can find a
clique such that the projection to the corresponding module is
non-zero, then we know that all characteristic vectors for maximum
independent sets are in the direct sum of the trivial module and the
standard module. Finally, we need to check that the matrix $M$ for $G$
has full rank. It is does, then the group has the strict EKR property.

This method has been applied to the symmetric group, alternating group
and $\PGL(2,q)$. In the next section we apply it to the Mathieu groups
and all 2-transitive groups on 20 or fewer points.

\section{EKR for the Mathieu groups}\label{sporadic}

In this section, using the module method, we establish the strict EKR
property for the $2$-transitive Mathieu groups. Since the family of Mathieu groups is
finite, the main approach of this problem uses a computer program to
show that the conditions of Theorem~\ref{module_method_thm} hold. All
of these programs have been run in the \textbf{GAP} programming system
\cite{GAP4}.

Following the standard notation, we will denote the Mathieu group of
degree $n$ by $M_n$. Note that $M_n\leq \sym(n)$ and we consider the
natural action of $M_n$ on the set $[n]$, as
usual. Table~\ref{transitivity_mathieu} lists some of the properties
of the Mathieu groups which will be useful for our purpose (see
\cite{cameron1999permutation} for more details).

\begin{table}[ht!]
\[
\begin{tabular}{|c|c|c|c|}
\hline
{\bf Group} & {\bf Order} & {\bf Transitivity} & {\bf Simplicity}\\
\hline
$M_{10}$ & 720 & sharply 3-transitive & not simple\\
\hline
$M_{11}$ & 7920 & sharply 4-transitive & simple\\
\hline
$M_{12}$ & 95040 & sharply 5-transitive & simple\\
\hline
$M_{21}$ & 20160 & 2-transitive & simple\\
\hline
$M_{22}$ & 443520 & 3-transitive & simple\\
\hline
$M_{23}$ & 10200960 & 4-transitive & simple\\
\hline
$M_{24}$ & 244823040 & 5-transitive & simple\\
\hline
\end{tabular}
\]
\caption{Order and transitivity table for Mathieu groups}
\label{transitivity_mathieu}
\end{table}

A simple computer program confirms the following.

\begin{lem}\label{standard_is_least_mathieu} Let $n\in
  \{10,11,12,21,22,23,24\}$. For each of the groups $M_n$ the least
  eigenvalue is $-\frac{\dd_{M_n}}{n-1}$ and the standard representation
  is the only representation with this eigenvalue.\qed
\end{lem}

This implies that condition (a) holds by the ratio bound. Furthermore,
condition (b) also holds for all the Mathieu groups, since the standard
module is the entire $-\frac{\dd_{M_n}}{n-1}$-eigenspace.

Finally, we need to confirm that condition (c) also holds for all the
Mathieu groups.  This requires showing that the matrix $M$ for each of
the Mathieu groups has full rank. For small $n$ we can do this with a
computer program.

\begin{lem}\label{M_10_11_12_21} 
If $G = M_n$ and $n\in \{10,11,12,21\}$, then the matrix $M$ for $G$ has full rank.\qed
\end{lem}

For $n = 22,23,24$ the Mathieu groups are too large to quickly check
the rank of the matrix $M$ using a computer, so instead we determine
the entries of the matrix $M^\top M$. The entries of this matrix can
be expressed as the linear combination of the identity matrix and the
adjacency matrix of a graph that we define next.

For $n >3$ define a graph $X_n$, which we call the \textsl{pairs
  graph}.  For any $n>3$, the vertices of $X_n$ are all the ordered
pairs $(i,j)$, where $i,j\in[n-1]$ and $i\neq j$; the vertices $(i,j)$
and $(k,\ell)$ are adjacent in $X_n$ if and only if either
$\{i,j\}\cap\{k,\ell\}=\emptyset$, ($i=\ell$ and $j\neq k$) or ($i\neq
\ell$ and $j=k$). The graph $X_n$ is regular of valency
$(n-2)(n-3)$. Note that the vertices of the pairs graph $X_n$ are the
pairs from $[n-1]$. The next lemma has been proved, using a version of
the ratio bound for cliques, in \cite{2013arXiv1302.7313A}.

\begin{lem}\label{least_eval_of_X_n}
 For any $n>3$, the least eigenvalue of the pairs graph $X_n$ is at least $-(n-3)$.\qed
\end{lem}

\begin{lem}\label{M_22} 
The matrix $M$ for the group $M_{22}$ has full rank.
\end{lem}
\proof Let $\mathcal{C}_{22}$ be one of
  the (two) conjugacy classes of $M_{22}$ whose elements are product
  of two disjoint $11$-cycles. Set $N=M_{\mathcal{C}_{22}}^\top
  M_{\mathcal{C}_{22}}$. Using a computer code we can establish
\[
N=1920\,I\,\,+\,\,96\,A(X_n).
\]
Lemma~\ref{least_eval_of_X_n} shows that the least eigenvalue of
$N$ is at least $1920-96(19)=96$. This shows that $N$ is non-singular
and we are done.\qed

\begin{lem}\label{M_23_24}
 The matrix $M$ for the group $M_n$ has full rank for $n\in \{23,24\}$.
\end{lem}
\proof
  Let $\cc_{23}$ be one of the two conjugacy classes of $M_{23}$ of
  permutations that are $23$-cycles and let $\cc_{24}$ be the only
  conjugacy class of $M_{24}$ whose elements are the product of two
  disjoint $12$-cycles. Set $t_n=|\cc_n|$, for $n=23,24$. Assume
  $M_{\cc_n}$ to be the submatrix of $M$, with the rows labeled by
  $\cc_n$ and set $N_n=M_{\cc_n}^\top M_{\cc_n}$, for $n=23,24$.  We
  now calculate the entries of $N$. Since $M_n$ is $4$-transitive, the
  entry $((a,b),(c,d))$ in $N_n$ depends only on the intersection of
  $\{a,b\}$ and $\{c,d\}$. To see this, consider the pairs
  $(a,b),(c,d)$ from $[n-1]$.  If an element $\pi\in \cc_n$ maps
  $a\mapsto b$ and $c\mapsto d$, then for any pairs $(a',b'),(c',d')$
  of elements of $[n-1]$, the permutation $g^{-1}\pi g\in \cc_n$ maps
  $a'\mapsto b'$ and $c'\mapsto d'$, where $g\in M_n$ is a permutation
  which maps $(a',b',c',d')$ to $(a,b,c,d)$. Therefore, we have
\begin{equation}\label{entries_of_N}
\hspace{-1cm}(N_n)_{(a,b),(c,d)}=
\begin{cases}
(N_n)_{(1,2),(1,2)}, &\quad \text{if } (c,d)=(a,b);\\
(N_n)_{(1,2),(2,1)}, & \quad\text{if } (c,d)=(b,a);\\
(N_n)_{(1,2),(2,3)}, & \quad\text{if }  a\neq d \quad\text{and}\quad b=c;\\
(N_n)_{(1,2),(2,3)}, & \quad\text{if }  a=d \quad\text{and}\quad b\neq c;\\
(N_n)_{(1,2),(3,4)}, & \quad\text{if } a,b,c,d\quad\text{are distinct}.\\
\end{cases}
\end{equation}
Because of the $2$-transitivity of $M_n$, we have
$(N_n)_{(1,2),(1,2)}=\frac{t_n}{n-1}$ and using a simple computer code we
can check that
\[(N_n)_{(1,2),(2,3)}=(N_n)_{(1,2),(3,4)}=\frac{t_n}{(n-1)(n-2)}.\]
 Also since elements of $\cc_n$ do not include a cycle of length $2$ in their cycle decomposition, we have $(N_n)_{(1,2),(2,1)}=0$. Thus we can re-write (\ref{entries_of_N}) as
\[
(N_n)_{(a,b),(c,d)}=
\begin{cases}
\frac{t_n}{n-1}, & \quad\text{if } (c,d)=(a,b);\\
0, & \quad\text{if } (c,d)=(b,a);\\
\frac{t_n}{(n-1)(n-2)}, &\quad \text{otherwise }.\\
\end{cases}
\]
This means that one can write
\[
N_n=\frac{t_n}{n-1}\,I\,\,+\,\,\frac{t_n}{(n-1)(n-2)}\,A(X_n),
\]
where $A(X_n)$ is the adjacency matrix of the pairs graph $X_n$. Then Lemma~\ref{least_eval_of_X_n} shows that the least eigenvalue of $N_n$ is at least 
\[
\frac{t_n}{n-1}\left(1-\frac{n-3}{n-2}\right)>0.
\]
We conclude that $N_n$ and, consequently, the matrix $M$ are full rank.\qed

Putting these together we have the that all the Mathieu groups has the strict EKR property.

\begin{thm}\label{main_mathieu}
The Mathieu groups $M_n$, for $n\in \{10,11,12,21,22,23,24\}$, have the strict EKR property.
\end{thm}

\section{2-transitive groups with small degree} 

In \ref{appA}, we present Table~\ref{small_perm_groups} in which we
record the results of applying the module method to the $2$-transitive
permutation groups of degree at most $20$.  For each of these groups
we first calculate the eigenvalues of the derangement graph. If the
standard character gives the least eigenvalue, then we know by
Theorem~\ref{ratio2} that the group has the EKR property.
If it is the only character that gives this eigenvalue,
then condition (b) of \ref{module_method_thm} holds as well.

If the standard character does not produce the least eigenvalue, then
we then check for sharply-transitive subgroups. If one exists, then
according to Theorem~\ref{clique_coclique_bound}, the group has the
EKR property (we did not check for sharply transitive sets, so it is
still possible that the group may have a clique of size $n$). If we do
find such cliques, we then check the projections of these cliques to
the different modules. If for every module, other than the trivial and
the standard, we can find a clique of size $n$ whose projection to the
module is non-zero, then condition (b) of \ref{module_method_thm}
holds.

Finally, for all the groups we check if the matrix $M$ for $G$ has full
rank (condition (c)). All these steps are implemented by a
\textbf{GAP} program.  Below we discuss two special cases.

If the eigenvalues of $\Gamma_G$ are $\dd_G$ and $-1$, then $\Gamma_G$
is the union of complete graphs. This means that the eigenvalue on the
standard module is $-1$ or that
\[
-\frac{\dd_G}{n-1} = -1.
\]
From this we can conclude that the degree of $\Gamma_G$ is $n-1$ and
that it is the union of $\frac{|G|}{n}$ copies of the complete graph
on $n$ vertices. This means that $G$ has the EKR property, but,
provided that $n >3$, it does not have the strict EKR property. To see
this, simply note that this graph has $n^{\frac{|G|}{n}}$ independent
sets of size $\frac{|G|}{n}$.  Since $G$ is $2$-transitive and,
provided that $n>3$, we have that $\frac{|G|}{n} >2$ and there are
more maximum independent sets than cosets of a point-stabilizer.  If
$n=3$ then the only group to consider is $\sym(3)$, which can easily
been seen to have the strict EKR property (this case is so small that it is
possible to list all the independent sets).

In \cite{MeagherS11} is it shown that the group $\PGL(n,q)$ has the
EKR property. This is due to the fact that these groups all contain a
Singer cycle; these cycles form large cliques in $\Gamma_{\PGL(n,q)}$
and the result follows from the clique-coclique bound.  But, provided
that $n >2$, these groups never have the strict EKR property. This
holds since the stabilizer of a hyperplane is an intersecting set with
size equal to the stabilizer of a point.

\section{Further Work}\label{future}

In the table in the appendix we have that every $2$-transitive group
with degree 20 or less has the EKR property. We would like to
determine if every $2$-transitive group has the EKR property.

We would also like to apply this method to other families of
$2$-transitive groups. For example, the group $\PSL(2,q)$ has the EKR
property and condition (b) holds (see~\cite{BahmanAhmadi} for
details). Furthermore, a computer program has shown that for every $q
\leq 30$, the matrix $M$ has full rank. But, we have been unable to
prove that $M$ has full rank for every $q$.

The next step will be to test if the groups for which the matrix $M$
is not full rank have the strict EKR property. One approach is to
search for the maximum independent sets in the derangement graphs. The
problem of finding maximum independent sets is NP-hard so this would
take a long time. Another approach would be to look for some
``natural'' independent sets.  For example, the stabilizer of a set of
points may form a maximum independent set, or perhaps the set of all
permutations that fix at least two points from a set of three points
would be the largest independent set.

The next open problem that we plan to work on is to determine if there
is any structure to the maximum intersecting sets in groups that do
not have the strict EKR property. For example, it has been
conjectured~\cite{MeagherS11} that the maximum independent set in $\Gamma_{\PGL(3,q)}$
are either stabilizers of a point or the stabilizer of a hyperplane.
It is not clear if the module method will be useful for this problem.

Our results for the subgroups with degree $16$ are not as satisfying. For many of
these groups the matrix $M$ does not have full rank, and we suspect
that the strict EKR property does not hold. It would be interesting to
look for non-canonical independent sets of maximum size in the
derangement graphs of these groups. Also, we would like to know if it is
possible to determine the structure of the derangement graph of the
group $G \rtimes \zz_2$ when we know that structure of the derangement
graph of $G$.

Finally, we plan to look more closely at the groups for which the
matrix $M$ has full rank, but we have not yet been able to prove that
the characteristic vector of every maximum independent set is in the
sum of the trivial and standard module. In fact, for all of the
$2$-transitve that groups we have considered, the maximum independent
sets are either all in the standard module or the derangement graph is
the union of complete graphs. It would be interesting to determine if
there are groups for which this does not hold.

\appendix

\section{Module Method for Small Groups}
\label{appA}

In this appendix we present a table of our results from applying the
the module method to all $2$-transitive groups with
degree at most $20$.

This work was implemented by a program in
\textbf{GAP}. Note that since all the groups $\sym(n)$ and $\alt(n)$
have the strict EKR property, they are excluded in the table.  In the
table we use the following terminology:

\begin{itemize}
\item \textbf{\textit{n}}: degree of the group;
\item {\bf least}: a ``Yes'' in this column means that the least
  eigenvalue of the derangement graph is given by the standard
  character;
\item {\bf $n$-clique}: a ``Yes'' in this column means that our program
  has found a clique of size $n$ in $\Gamma_n$ (hence the
  clique-coclique bound holds with equality); the symbol ``--'' means
  that we don't try to find a maximum clique, and the symbol ``?''
  means that the program failed to find such a clique (but not that one does not exist!);
\item {\bf EKR}: a ``Yes'' in this column means that the group has the
  EKR property, i.e. condition (a) of the module method holds;
\item {\bf unique}: a ``Yes'' in this column means that the standard
  character is the only character giving the least eigenvalue; hence
  condition (b) of the module method holds; An ``N/A'' means that this
  condition is not applicable since the standard character does not
  give the least eigenvalue.
\item {\bf module by clique}: a ``Yes'' in this column means that
  using cliques of size $n$ and Corollary~\ref{at_most_one_non-zero}
  we know that the characteristic vector of any maximum independent
  set of $\Gamma_G$ lies in the direct sum of the trivial and the
  standard characters of $G$; hence condition (b) of the module method
  holds; the symbol ``--'' means that we don't verify this, and the
  symbol ``?''  means that the program could not find suitable cliques
  to prove that condition (b) holds;
\item {\bf rank}: a ``Yes'' in this column means that the matrix $M$
  for the group $G$ has full rank, i.e. condition (c) of the module
  method holds; 
\item {\bf strict}: a ``Yes'' in this column means that $G$ has the
  strict EKR property; the symbol ``?''  means that the program could
  not verify this. In all the cases where we have a ``No'', either the
  derangement graph is the union of complete graphs or the group is
  $\PGL(3,q)$ for some $q$;
\end{itemize}

\newpage
\begin{landscape}
\begin{center}
\small
\begin{longtable}{|c|c|c|c|c|c|c|c|c|c|c|}
\caption[EKR and strict EKR property for small $2$-transitive groups]{EKR and strict EKR property for small $2$-transitive groups} \label{small_perm_groups} \\
\hline
 \multicolumn{1}{|c|}{$n$} & \multicolumn{1}{|c|}{Group} & \multicolumn{1}{|c|}{size} & \multicolumn{1}{|c|}{least}& \multicolumn{1}{|c|}{$n$-clique}& \multicolumn{1}{|c|}{{\bf EKR}}& \multicolumn{1}{|c|}{unique}& \multicolumn{1}{|c|}{module by clique} & \multicolumn{1}{|c|}{rank}& \multicolumn{1}{|c|}{{\bf strict}}\\ \hline\hline 
\endfirsthead

\multicolumn{10}{c}%
{{\bfseries \tablename\ \thetable{} -- continued from previous page}} \\
\hline 
\multicolumn{1}{|c|}{$n$} & \multicolumn{1}{|c|}{Group} & \multicolumn{1}{|c|}{size} & \multicolumn{1}{|c|}{least}& \multicolumn{1}{|c|}{$n$-clique}& \multicolumn{1}{|c|}{{\bf EKR}}& \multicolumn{1}{|c|}{unique}& \multicolumn{1}{|c|}{module by clique} & \multicolumn{1}{|c|}{rank}& \multicolumn{1}{|c|}{{\bf strict}}\\ \hline \hline
\endhead

\hline \multicolumn{10}{|r|}{{Continued on next page}} \\ \hline
\endfoot
\hline \hline
\endlastfoot

$5$ & $\zz_5\rtimes \zz_4$ & $20$  & Yes & -- & Yes & Yes & -- & No & No\\
\hline
$6$ & $\PGL(2,5)$ & $120$  & Yes & Yes & Yes & No & -- & Yes & Yes\\
\hline
$6$ & $\alt(5)$ & $60$  & Yes & -- & Yes & Yes & -- & Yes & Yes\\
\hline
$7$ & $\PGL(3,2)$ & $168$  & Yes & Yes & Yes & No & ? & No & No\\
\hline
$7$ & $(\zz_7\rtimes \zz_3)\rtimes \zz_2$ & $42$  & Yes & -- & Yes & Yes & -- & No & No\\
\hline
$8$ & $(\zz_2\times\zz_2\times\zz_2)\rtimes \PSL(3,2)$ & $1344$  & Yes & -- & Yes & Yes & -- & Yes & Yes\\
\hline
$8$ & $\PGL(2,7)$ & $336$  & Yes & Yes & Yes & No & ? & Yes & Yes \\
\hline
$8$ & $((\zz_2\times\zz_2\times\zz_2)\rtimes \zz_7)\rtimes\zz_3$ & $168$  & No & Yes & Yes & N/A & Yes & Yes & Yes\\
\hline
$8$ & $\PSL(3,2)$ & $168$  & Yes & -- & Yes & Yes & -- & Yes & Yes\\
\hline
$8$ & $(\zz_2\times\zz_2\times\zz_2)\rtimes \zz_7$ & $56$  & Yes & -- & Yes & Yes & -- & No & No\\
\hline 
$9$ & $\PSL(2,8)\rtimes\zz_3$ & $1512$  & Yes & -- & Yes & Yes & -- & Yes & Yes\\
\hline
$9$ & $(((\zz_3\times\zz_3)\rtimes Q_8)\rtimes\zz_3)\rtimes \zz_2$ & $432$  & Yes & Yes & Yes & No & ? & Yes & ?\\
\hline
$9$ & $((\zz_3\times\zz_3)\rtimes Q_8)\rtimes \zz_3$ & $216$  & No & Yes & Yes & N/A & ? & No & ?\\
\hline 
$9$ & $\PSL(2,8)$ & $504$  & Yes & -- & Yes & Yes & -- & Yes & Yes\\
\hline
$9$ & $((\zz_3\times\zz_3)\rtimes \zz_8)\rtimes \zz_2$ & $144$  & No & Yes & Yes & N/A & ? & No & ?\\
\hline
$9$ & $(\zz_3\times\zz_3)\rtimes \zz_8$ & $72$ & Yes & -- & Yes & Yes & -- & No & No\\
\hline
$9$ & $(\zz_3\times\zz_3)\rtimes Q_8$ & $72$ & Yes & -- & Yes & Yes & -- & No & No\\
\hline
$10$ & $(\alt(6)\times\zz_2)\rtimes \zz_2$ & $1440$ & No & Yes & Yes & N/A & ?& Yes & ?\\
\hline
$10$ & $M_{10}$ & $720$  & Yes & -- & Yes & Yes & -- & Yes & Yes\\
\hline
$10$ & $\alt(6)\cdot \zz_2$ & $720$  & Yes & ? & Yes & No & ? & Yes & ?\\
\hline
$10$ & $\PGL(2,9)$ & $720$ & Yes & Yes & Yes & No & ?& Yes & Yes\\
\hline
$10$ & $\alt(6)$ & $360$  & Yes & -- & Yes & Yes & -- & Yes & Yes\\
\hline
$11$ & $M_{11}$ & $7920$  & Yes & -- & Yes & Yes & -- & Yes & Yes\\
\hline
$11$ & $\PSL(2,11)$ & $660$  & Yes & -- & Yes & Yes & -- & No & ?\\
\hline
$11$ & $(\zz_{11}\rtimes \zz_5)\rtimes \zz_2$ & $110$  & Yes & -- & Yes & Yes & -- & No & No\\
\hline
$12$ & $M_{12}$ & $95040$  & Yes & -- & Yes & Yes & -- & Yes & Yes\\
\hline
$12$ & $M_{11}$ & $7920$  & Yes & -- & Yes & Yes & -- & Yes & Yes\\
\hline
$12$ & $\PGL(2,11)$ & $1320$  & Yes & Yes & Yes & No & -- & Yes & Yes \\
\hline
$12$ & $\PGL(2,11)$ & $660$  & Yes & -- & Yes & Yes & -- & Yes & Yes\\
\hline
$13$ & $\PSL(3,3)$ & $5616$  & Yes & -- & Yes & Yes & -- & No & No\\
\hline
$13$ & $(\zz_{13}\rtimes\zz_4)\rtimes \zz_3$ & $156$  & Yes & -- & Yes & Yes & -- & No & No\\
\hline
$14$ & $\PGL(2,13)$ & $2184$  & Yes & Yes & Yes & No & -- & Yes & Yes \\
\hline
$14$ & $\PSL(2,13)$ & $1092$  & Yes & -- & Yes & Yes & -- & Yes & Yes\\
\hline
$15$ & $\alt(8)$ & $20160$  & Yes & -- & Yes & Yes & Yes & No & ?\\
\hline
$15$ & $\alt(7)$ & $2520$  & Yes & -- & Yes & Yes & -- & No & ?\\
\hline
$16$ & $(\zz_2\times\zz_2\times\zz_2\times\zz_2)\rtimes \alt(8)$ & $322560$  & Yes & Yes & Yes & Yes & -- & Yes & Yes\\
\hline
$16$ & $((\zz_2\times\zz_2\times\zz_2\times\zz_2)\rtimes \alt(6))\rtimes\zz_2$ & $11520$  & No & Yes & Yes & N/A & Yes & No & ?\\
\hline
$16$ & $(((\zz_2\times\zz_2\times\zz_2\times\zz_2)\rtimes \alt(5))\rtimes\zz_3)\rtimes\zz_2$ & $5760$  & No & Yes & Yes & N/A & Yes & Yes & Yes\\
\hline
$16$ & $((\zz_2\times\zz_2\times\zz_2\times\zz_2)\rtimes \alt(5))\rtimes\zz_3$ & $2880$  & Yes & Yes & Yes & No & ? & Yes & ?\\ 
\hline
$16$ & $(\zz_2\times\zz_2\times\zz_2\times\zz_2)\rtimes \alt(7)$ & $40320$  & Yes & Yes & Yes & Yes & Yes & Yes & Yes\\
\hline
$16$ & $(\zz_2\times\zz_2\times\zz_2\times\zz_2)\rtimes \alt(6)$ & $5760$  & Yes & Yes & Yes & No & Yes & No & ?\\
\hline
$16$ & $((\zz_2\times\zz_2\times\zz_2\times\zz_2)\rtimes \alt(5))\rtimes \zz_2$ & $1920$  & No & Yes & Yes & N/A & Yes & No & ?\\
\hline
$16$ & $(\zz_2\times\zz_2\times\zz_2\times\zz_2)\rtimes \alt(5)$ & $960$  & No & Yes & Yes & N/A & Yes & No & ?\\
\hline
$16$ & $(((\zz_2\times\zz_2\times\zz_2\times\zz_2)\rtimes \zz_5)\rtimes\zz_3)\rtimes\zz_4$ & $960$  & No & Yes & Yes & N/A & Yes & No & ?\\
\hline
$16$ & $(((\zz_2\times\zz_2\times\zz_2\times\zz_2)\rtimes \zz_5)\rtimes\zz_3)\rtimes\zz_2$ & $480$  & No & Yes & Yes & N/A & ? & No & ?\\
\hline
$16$ & $((\zz_2\times\zz_2\times\zz_2\times\zz_2)\rtimes \zz_5)\rtimes\zz_3$ & $240$  & Yes & -- & Yes & Yes & Yes & No & No\\
\hline
$17$ & $\PSL(2,16)\rtimes \zz_4$ & $16320$  & Yes & -- & Yes & Yes & -- & Yes & Yes\\
\hline
$17$ & $\PGL(2,16)$ & $8160$  & Yes & -- & Yes & Yes & -- & Yes & Yes\\
\hline
$17$ & $\PSL(2,16)$ & $4080$  & Yes & -- & Yes & Yes & -- & Yes & Yes\\
\hline
$17$ & $\zz_{17}\rtimes\zz_{16}$ & $272$  & Yes & -- & Yes & Yes & -- & No & No\\
\hline
$18$ & $\PGL(2,17)$ & $4896$  & Yes & -- & Yes & No & -- & Yes & Yes \\
\hline
$18$ & $\PSL(2,17)$ & $2448$  & Yes & -- & Yes & Yes & -- & Yes & Yes\\
\hline
$19$ & $(\zz_{19}\rtimes\zz_{9})\rtimes\zz_2$ & $342$  & Yes & -- & Yes & Yes & -- & No & No\\
\hline
$20$ & $\PGL(2,19)$ & $6840$  & Yes & -- & Yes & No & -- & Yes & Yes \\
\hline
$20$ & $\PSL(2,19)$ & $3420$  & Yes & -- & Yes & Yes & -- & Yes & Yes\\
\hline
\end{longtable}
\end{center}
\end{landscape}


\end{document}